\documentclass[11pt]{article}  

\usepackage[dvips]{epsfig}
\usepackage{latexsym}
\usepackage{amsfonts,amsmath,amssymb}
\newtheorem{corollary}{Corollary}
\newtheorem{lemma}{Lemma}
\newtheorem{example}{Example}
\newtheorem{theorem}{Theorem}
\newtheorem{proposition}{Proposition}

\newtheorem{conjecture}{Conjecture}
\newtheorem{definition}{Definition}
\newcommand{\qed}{\mbox{$\Diamond$}\vspace{\baselineskip}}
\newenvironment{proof}{\noindent{\bf Proof:}}{\qed}
\newcommand{\Var}{\hbox {Var}}

\usepackage{lineno}

\begin{document}

\title{On Three Different Notions of Monotone Subsequences} 

\author{Mikl\'os B\'ona \\
Department of Mathematics \\
University of Florida\\
 Gainesville FL 32611-8105\\
bona@math.ufl.edu \thanks{Partially supported
by an NSA Young Investigator Award.}}

\date{}
\maketitle

\begin{abstract}
We review how the monotone pattern compares to other patterns in terms
of enumerative results on pattern avoiding permutations. We consider three
natural definitions of pattern avoidance, give an overview of classic and
recent formulas, and provide some new results related to limiting 
distributions.
\end{abstract} 

\section{Introduction}
Monotone subsequences in a permutation $p=p_1p_2\cdots p_n$ has been the 
subject of vigorous research for over sixty years. In this paper, we will
review three different lines of work. In all of them, we will consider
increasing subsequences of a permutation of length $n$ that have a 
{\em fixed} length $k$. This is in contrast to another line of work,
started by Ulam more than sixty years ago, in which the distribution of
the {\em longest} increasing subsequence of a random permutation has been
studied. That direction of research has recently reached a high point in
the article \cite{Deift} of Baik, Deift and Johansson. 

The three directions we consider are distinguished by their definition
of monotone subsequences. We can simply require that $k$ entries of a 
permutation increase from left to right, or we can in addition require that 
these $k$ entries be in consecutive positions, or we can even require that
in  they be consecutive integers {\em and} be in consecutive positions.

\section{Monotone Subsequences with No Restrictions}
The classic definition of pattern avoidance for permutations is as follows.
Let $p=p_1p_2\cdots p_n$ be a permutation, let $k<n$, and let $q=q_1q_2\cdots
q_k$ be another permutation. We say that $p$ {\em contains} $q$ as a pattern
if there exists a subsequence $1\leq i_1<i_2<\cdots <i_k\leq n$
 so that for all indices
$j$ and $r$, the inequality $q_j<q_r$ holds if and only if the inequality
$p_{i_j}<p_{i_r}$ holds. If $p$ does not contain $q$, then we say
that $p$ {\em avoids} $q$.  In other words, $p$ contains $q$ if $p$ has
a subsequence of entries, not necessarily in consecutive positions, which
relate to each other the same way as the entries of $q$ do. 

\begin{example} The permutation 3174625 contains the pattern 123. Indeed,
consider the first, fourth, and seventh entries.
\end{example}

In particular, $p$ contains the monotone pattern $\alpha_k=12\cdots k$ if
and only if $p$ contains an increasing subsequence of length $k$. The
elements of this increasing subsequence do not have to be in consecutive
positions. 

The enumeration of permutations avoiding a given pattern is a fascinating
subject. Let $S_n(q)$ denote the number of permutations of length $n$ 
(or, in what follows, $n$-permutations) that
avoid the pattern $q$. 

\subsection{Patterns of Length Three}
Among patterns of length three, there is no difference
between the  monotone pattern and other patterns as far as $S_n(q)$ is
concerned. This is the content of our first theorem.

\begin{theorem} \label{allsix}
Let $q$ be any pattern of length three, and let $n$ be any positive integer.
Then $S_n(q)=C_n={2n\choose n}/(n+1)$. In other words, $S_n(q)$ is the $n$th
Catalan number.
\end{theorem}

\begin{proof} If $p$ avoids $q$, then the reverse of $p$ avoids the reverse
of $q$, and the complement of $p$ avoids the complement of $q$. Therefore,
$S_n(123)=S_n(321)$ and $S_n(132)=S_n(231)=S_n(213)=S_n(312)$. 

The fact that $S_n(132)=S_n(123)$ is proved using the well-known 
Simion-Schmidt bijection \cite{Simion}. In a permutation,
let us call an entry a {\em left-to-right} minimum if it is smaller than
every entry on its left. For instance, the left-to-right minima of 
4537612 are the entries 4, 3, and 1. 
 
 Take an $n$-permutation $p$ of length $n$ that
avoids 132, keep its left-to-right minima fixed, and arrange all other entries
in decreasing order in the positions that do not belong to left-to-right
minima, to get the permutation $f(p)$. For instance, if $p=34125$, then
$f(p)=35142$. Then $f(p)$ is a union of two decreasing sequences, so it is
123-avoiding. Furthermore, $f$ is a bijection between the two relevant set
of permutations. Indeed, if $r$ is a permutation counted by $S_n(123)$, then
$f^{-1}(r)$ is obtained by keeping the left-to-right minima of $r$ fixed, 
and rearranging the remaining entries so that moving from left to right, 
each slot is filled by the smallest remaining entry that is larger than the
closest left-to-right minimum on the left of that position.

In order to prove that $S_n(132)=C_n$, just note that in a 132-avoiding
$n$-permutation, any entry to the left of $n$ must be smaller than any entry
to the right of $n$. Therefore, if $n$ is in the $i$th position, then there
are $S_{i-1}(132)S_{n-i}(123)$ permutations of length $n$ that avoid 132.
Summing over all $i$, we get the recurrence
\[S_{n}(132)= \sum_{i=0}^{n-1}S_{i-1}(132)S_{n-i}(132),\]
which is the well-known recurrence for Catalan numbers. 
\end{proof}

\subsection{Patterns of Length Four}
When we move to longer patterns, the situation becomes much more complicated
and less well understood. In his doctoral thesis \cite{West},
 Julian West published the 
following numerical evidence.

\begin{itemize}
\item for $S_n(1342)$, and $n=1,2,\cdots ,8$, 
 we have 1, 2, 6, 23, 103, 512, 2740, 15485
\item for $S_n(1234)$, and $n=1,2,\cdots ,8$, 
 we have  1, 2, 6, 23, 103, 513, 2761, 15767
\item for $S_n(1324)$, and $n=1,2,\cdots ,8$, 
 we have  1, 2, 6, 23, 103, 513, 2762, 15793.
\end{itemize}

These data are startling for at least two reasons. First, the numbers
$S_n(q)$ are no longer independent of  $q$; there
are some patterns of length four that are easier to avoid than others.
Second, the monotone pattern 1234, special as it is, 
 does not provide the minimum or the
maximum value for $S_n(q)$. We point out that for 
each $q$ of the other 21 patterns of length four, it is known that the
sequence $S_n(q)$ is identical to one of the three sequences $S_n(1342)$,
$S_n(1234)$, and  $S_n(1324)$. See \cite{bona}, Chapter 4, for more
details.  

Exact formulas are known for two of the above three sequences. 
For the monotone pattern, Ira Gessel gave a formula using symmetric functions.

\begin{theorem} \cite{GesselF}, \cite{GesselP}
For all positive integers $n$, the identity
\begin{eqnarray} S_n(1234) & = & 
2\cdot\sum_{k=0}^n {2k\choose k}{n\choose k}^2 
\frac{3k^2+2k+1-n-2nk}{(k+1)^2 (k+2) (n-k+1) } \\
& = & \frac{1}{(n+1)^2(n+2)} \sum_{k=0}^n {2k\choose k}{n+1\choose k+1}
{n+2\choose k+1}.\end{eqnarray}
\end{theorem}

The formula for $S_n(1342)$ is due to the present author \cite{Bona1342},
 and is quite surprising. 
\begin{theorem} \label{exactf} For all positive integers $n$, we have
\begin{eqnarray*} S_n(1342) & = & (-1)^{n-1} \cdot \frac{(7n^2-3n-2)}{2} \\
 & + & 3\sum_{i=2}^n (-1)^{n-i} \cdot
 2^{i+1}\cdot \frac{(2i-4)!}{i!(i-2)!}\cdot 
{{n-i+2\choose 2}}
. \end{eqnarray*}
\end{theorem}

This result is unexpected for two reasons. First, it shows that 
$S_n(1342)$ is not simply less than $S_n(1234)$ for every $n\geq 6$; it is 
{\em much less}, in a sense that we will explain in Subsection \ref{swlimits}.
For now, we simply state that while $S_n(1234)$ is ``roughly'' $9^n$, 
the value of $S_n(1342)$ is``roughly'' $8^n$. Second, the formula is, in
some sense, simpler than that for $S_n(1234)$. Indeed, it follows 
from Theorem \ref{exactf} that the ordinary generating
function of the sequence $S_n(1342)$ is 
\[H(x)=\sum_{i\geq 0}F^i(x)=\frac{1}{1-F(x)}=
\frac{32x}{-8x^2+20x+1-(1-8x)^{3/2}}.\]
This is an {\em algebraic} power series. On the other hand,
it is  known (Problem Plus 5.10 in \cite{bona} that the 
ordinary generating
function of the sequence $S_n(1234)$ is {\em not} algebraic.
So permutations avoiding the 
 monotone pattern are not even the {\em nicest} among permutations avoiding
a given pattern, in terms of the generating functions that count them. 

There is no known formula for the third sequence, that of the numbers
$S_n(1324)$. However, the following inequality is known \cite{bonathesis}.
 
\begin{theorem} For all integers $n\geq 7$, the inequality
\[S_n(1234) < S_n(1324) \]
holds. 
\end{theorem}

\begin{proof} Let us call an entry of a permutation a {\em 
right-to-left maximum} if it is larger than all entries on its right.
So Let us say that two $n$-permutations are in the same class 
if they have the same left-to-right minima, and they are in the same 
positions, and they have the same right-to-left maxima, and they are
in the same positions as well. 
For example, $51234$ and $51324$ are in the same
class, but $z=24315$  and
$v=24135$ are not, as the third entry of $z$ is not a
left-to-right minimum, whereas that of $v$ is.

It is straightforward to see that each non-empty class contains exactly one 
1234-avoiding permutation, the one in which the subsequence of entries that
are neither left-to-right minima nor right-to-left maxima is decreasing. 

It is less obvious that each class contains {\em at least one} 1324-avoiding
permutation. Note that if a permutation contains a 1324-pattern,
then we can choose such a pattern so that its first element is a left-to-right
 minimum and its last element is a right-to-left maximum. Take a 1324-avoiding
permutation, and take one of its 1324-patterns of the kind described in the
previous sentence. Interchange its second and third element. Observe that
this will keep the permutation within its original class. Repeat this
 procedure as long as possible. The procedure will stop after a finite
number of steps since each step decreases the number of inversions of the
permutation. When the procedure stops, the permutation at hand avoids 1324.

This shows that $ S_n(1234) \leq S_n(1324)$ for all $n$. If $n\geq 7$, then
the equality cannot hold since there is at least one class that contains more
than one 1324-avoiding permutation. For $n=7$, this is the class
$3*1*7*5$, which contains 3612745 and 3416725. For larger $n$, this class
can be prepended by $n(n-1)\cdots 8$ to get a suitable class. 
\end{proof}

It turns out again that $S_n(1324)$ is {\em much} larger than $S_n(1234)$. 
We will give the details in Subsection \ref{swlimits}.

\subsection{Patterns of Any Length}
For general $k$, there are some good estimates known 
for the value of $S_n(\alpha_k)$. The first one can be proved by
an elementary method.

\begin{theorem} \label{1rank} 
For all positive integers $n$ and $k>2$, we have
\[S_n(123\cdots k)\leq (k-1)^{2n}.\]
 \end{theorem}

\begin{proof}
Let us say that an entry $x$ of a permutation is of rank $i$ if it is
the end of an increasing subsequence of length $i$, but there is no increasing
subsequence of length $i+1$ that ends in $x$. Then for all $i$,
elements of rank $i$ must form a decreasing subsequence. Therefore, a
$q$-avoiding permutation can be decomposed into the union of $k-1$
decreasing subsequences. Clearly, there are at most
$(k-1)^{n}$ ways to partition
our $n$ entries into $k-1$ blocks. Then we have to place these blocks of
entries
somewhere in our permutation.  There are at most  $(k-1)^{n}$ ways
to assign  each position of the permutation 
to one of these blocks, completing the
proof.
\end{proof}

 Indeed, 
Theorem \ref{1rank}  has a stronger version, obtained by Amitaj Regev
\cite{Regev}. It needs heavy analytic
machinery, and therefore  will not be proved here.
We mention the result, however, as it shows that no matter what $k$ is, 
 the constant  $(k-1)^2$ in Theorem \ref{1rank} cannot 
be replaced by a smaller number, so the elementary estimate
of Theorem \ref{1rank} is optimal in some strong sense. We remind the 
reader that functions
$f(n)$ and $g(n)$ are said to be {\em asymptotically equal} if
$\lim_{n\rightarrow \infty} \frac{f(n)}{g(n)}=1$. 

\begin{theorem} \label{regev} \cite{Regev}
\label{monoton} For all $n$, $S_n(1234\cdots k)$
asymptotically equals \[\lambda_k \frac{(k-1)^{2n}}{n^{(k^2-2k)/2}} .\]  Here
\[ \lambda_k=\gamma_k^2 
\int\!\!\!\!\!\int\limits_{x_1\,\geq\,x_2\,\geq\,\cdots\,\geq\,x_k}
\!\!\!\!\!\!\!\!\!\!\!\!\!\cdots\int
[D(x_1,x_2,\cdots ,x_k)\cdot
e^{-(k/2)x^2 }] ^2 dx_1 dx_2\cdots dx_k, \]
where  $D(x_1,x_2,\cdots ,x_k)=\Pi_{i<j} (x_i-x_j)$, and
$\gamma_k= (1/\sqrt{2\pi})^{k-1}\cdot k^{k^2/2}. $
 \end{theorem}

\subsection{Stanley-Wilf Limits} \label{swlimits}

The following
 celebrated result of Adam Marcus and G\'abor Tardos \cite{Marcus} shows
that in general, it is very difficult to avoid any given pattern $q$.

\begin{theorem} \cite{Marcus}
 For all patterns $q$, there exists a constant $c_q$ so that
\begin{equation} S_n(q)\leq c_q^n.\end{equation}
\end{theorem}

It this not difficult to show using Fekete's lemma that the sequence
$\left(S_n(q) \right)^{1/n}$ is monotone increasing. The previous
theorem shows that it is bounded from above, leading to the following.

\begin{corollary} \label{limits}
For all patterns $q$, the limit 
\[L(q)=\lim_{n\rightarrow \infty} \left(S_n(q) \right)^{1/n}
\]
exists.
\end{corollary}

The real number $L(q)$ is called the {\em Stanley-Wilf} limit, or 
{\em growth rate} of the pattern $q$. In this terminology, Theorem
\ref{regev} implies that $L(\alpha_k)=(k-1)^2$. In particular,
$L(1234)=9$, while Theorem \ref{exactf} implies that $L(1342)=8$. So 
it is not simply easier to avoid 1234  than 1342, it is {\em exponentially}
easier to do so. 

Numerical evidence suggests that in the multiset of $k!$ real numbers 
$S_n(q)$, the numbers $S_n(\alpha_k)$ are much closer to the maximum than
to the minimum. This led to the plausible conjecture that 
for any pattern $q$ of length $k$, the inequality $L(q)\leq (k-1)^2$
holds. This would mean that while there are patterns of length $k$
 that are easier to avoid
than $\alpha_k$, there are none that are much easier to avoid, in the sense
of Stanley-Wilf limits.
 However, this conjecture has been disproved by the following
result of Michael Albert and al.

\begin{theorem} \cite{albert}
The inequality $L(1324)\geq 11.35$ holds.
\end{theorem}

In other words, it is not simply harder to avoid 1234 than 1324, 
  it is {\em exponentially} harder to do so. 

\subsection{Asymptotic Normality} \label{asymptotics}
In this section we change direction and  prove that the distribution of
the number of
copies of $\alpha_k$ in a randomly selected $n$-permutation converges
in distribution to a normal distribution. (For the rest of this
paper, when we say random permutation of length $n$, we always assume 
that each $n$-permutation is selected with probability $1/n!$.)
 Note that in the special case
of $k=2$, this is equivalent to 
the classic result that the distribution of inversions
in random permutations is asymptotically normal. See \cite{fulman} and its
references for various proofs of that result, or \cite{ngendes} for a
generalization.

We need to introduce some notation for transforms of the random variable
$Z$. Let $\bar{Z}=Z-E(Z)$, let $\tilde{Z}=\bar{Z}/\sqrt{\Var( Z)}$, and let
$Z_n\rightarrow N(0,1)$ mean that $Z_n$ converges in distribution to the 
standard normal variable. 

Our main tool in this section will be a theorem of Svante Janson
 \cite{janson}. In order to be able to state that theorem, we need the
following definition. 

\begin{definition}
Let $\{Y_{n,k}|k=1,2,\cdots ,N_n\}$ be an array of
 random variables.
 We say that a graph $G$ is 
a {\em dependency graph} for   $\{Y_{n,k}|k=1,2\cdots , N_n\}$
 if the following
two conditions are satisfied:
\begin{enumerate}
\item There exists a bijection between the random variables $Y_{n,k}$ and
the vertices of $G$, and
\item If $V_1$ and $V_2$ are two disjoint sets of vertices of $G$ so that
no edge of $G$ has one endpoint in $V_1$ and another one in $V_2$, then
the corresponding sets of random variables are independent.
\end{enumerate}
\end{definition}

Note that  the dependency graph of a
 family of variables is not unique. Indeed if $G$ is a dependency graph
for a family and $G$ is not a complete graph,
 then we can get other dependency graphs for the family
by simply adding new edges to $G$. 

Now we are in position to state Janson's theorem, the famous
{\em Janson dependency criterion}.

\begin{theorem} \cite{janson} \label{janson}
Let $Y_{n,k}$ be an array of random variables such that for all $n$, and
for all $k=1,2,\cdots ,N_n$, the inequality $|Y_{n,k}|\leq A_n$ holds for
some real number $A_n$, and that the maximum degree of a dependency
graph of $\{Y_{n,k} | k=1,2,\cdots ,N_n \}$ is $\Delta_n$. 

Set $Y_n=\sum_{k=1}^{N_n} Y_{n,k}$ and $\sigma_n^2= \Var ( Y_n)$. If there
is a natural number $m$ so that
\begin{equation} \label{jansencond}
N_n\Delta_n^{m-1} \left (\frac{A_n}{\sigma_n} \right )^m \rightarrow 0,
\end{equation}
as $n$ goes to infinity, then \[ \tilde{Y}_n \rightarrow N(0,1) .\]
\end{theorem}

Let us order the ${n\choose k}$ subwords of length $k$ of the permutation
$p_1p_2\cdots p_n$ linearly in some way. 
For $1\leq i\leq {n\choose k}$, let $X_{n,i}$ 
 be the indicator random
variable of the event that in a randomly selected permutation of length $n$,
the $i$th subword of length $k$ in the permutation $p=p_1p_2\cdots p_n$
is a $12\cdots k$-pattern. We will now verify that the family of the
$X_{n,i}$ satisfies all conditions of the Janson Dependency Criterion.

First, $|X_{n,i}|\leq 1$ for all $i$ and all $n$, since the $X_{n,i}$ are 
indicator random variables. So we can set $A_n=1$. Second, $N_n={n\choose k}$,
the total number of subwords of length $k$ in $p$. Third, if $a\neq b$, then
$X_a$ and $X_b$ are independent unless the corresponding subwords intersect.
For that, the $b$th subword must intersect the $a$th subword in $j$ entries, 
for some $1\leq j\leq k-1$. For a fixed $a$th subword, the number of 
ways that can happen is $\sum_{j=1}^{k-1} {k\choose j}{n-k\choose k-j}=
{n\choose k}-{n-k \choose k}-1$, where we used 
 the well-known Vandermonde identity to compute the sum.
Therefore, 
\begin{equation} \label{maxdegree} 
\Delta_n \leq {n\choose k}-{n-k \choose k}-1.
\end{equation}
In particular, note that (\ref{maxdegree}) provides an upper bound for
$\Delta_n$ in terms of a polynomial function of $n$ that is 
of degree $k-1$ since terms of degree
$k$ will cancel.

There remains the task of finding a lower bound for $\sigma_n$ that 
we can then use in applying Theorem \ref{janson}. Let $X_n=
\sum_{i=1}^{n\choose k} X_{n,i}$. We will show the following.

\begin{proposition} \label{varprop}
 There exists a positive constant $c$ so that
for all $n$, the inequality
\[\Var(X_n)\geq cn^{2k-1}\]
holds.
\end{proposition}

\begin{proof}
By linearity of expectation, we have
\begin{eqnarray} \label{variance}
\Var (X_n) & = & E(X_n^2) - (E(X_n))^2 \\
 & = & E \left (\left( \sum_{i=1}^{{n\choose k}} X_{n,i} \right )^2 \right )
-  \left (E \left (\sum_{i=1}^{{n\choose k}} X_{n,i} \right ) \right )^2 \\
 & = & E \left (\left( \sum_{i=1}^{{n\choose k}} X_{n,i} \right )^2 \right )
- \left( \sum_{i=1}^{{n\choose k}} E(X_{n,i}) \right )^2  \\
 \label{lastone} & = &  \sum_{i_1, i_2}
E(X_{n,i_1}X_{n,i_2})  - \sum_{i_1, i_2}
E(X_{n,i_1})E(X_{n,i_2}).
\end{eqnarray}

Let $I_1$ (resp. $I_2$) denote the $k$-element subword of $p$ indexed
by $i_1$, (resp. $i_2$). Clearly, it suffices to show that
\begin{equation} \label{simplified} \sum_{|I_1\cap I_2| \leq 1}
E(X_{n,i_1}X_{n,i_2}) - \sum_{i_1, i_2}
E(X_{n,i_1})E(X_{n,i_2}) \geq cn^{2k-1},\end{equation}
since the left-hand side of (\ref{simplified}) is obtained from the 
(\ref{lastone}) by removing the sum of some positive terms, that is,
the sum of all $E(X_{n,i_1}X_{n,i_2})$ where  $|I_1\cap I_2| >1$.

As $E(X_{n,i})=1/k!$ for each $i$, the sum with negative sign in 
(\ref{lastone}) is
\[ \sum_{i_1, i_2}
E(X_{n,i_1})E(X_{n,i_2}) ={n\choose k}^2 \cdot \frac{1}{k!^2},\]
which is a polynomial function 
 in $n$, of degree $2k$ and of leading coefficient
$\frac{1}{k!^4}$. As far as the summands in  (\ref{lastone}) with a positive
sign go, {\em most} of them are also equal to  $\frac{1}{k!^2}$. More
precisely, $E(X_{n,i_1}X_{n,i_2})=\frac{1}{k!^2}$ when 
 $I_1$ and $I_2$ are disjoint, and that happens for 
${n\choose k}{n-k\choose k}$ ordered pairs $(i_1,i_2)$
 of indices. The sum of these
summands is
\begin{equation}
\label{disjoint} d_n={n\choose k}{n-k\choose k} \frac{1}{k!^2},
\end{equation}
which is again a polynomial function in $n$, of degree $2k$ and with leading 
coefficient
$\frac{1}{k!^4}$. So summands
 of degree $2k$ will cancel out in (\ref{lastone}). (We will see in the next
paragraph that the summands we have not yet considered add up to a polynomial
of degree $2k-1$.)
In fact, considering the two types of summands we studied in
(\ref{lastone}) and (\ref{disjoint}), we see that they add up to 
\begin{eqnarray} 
{n\choose k}{n-k\choose k} \frac{1}{k!^2}-{n\choose k}^2  \frac{1}{k!^2}
& = & n^{2k-1} \frac{2{k\choose 2}-{2k-1\choose 2}}{k!^4}+O(n^{2k-2}) \\
\label{theeasy} & = & n^{2k-1} \frac{-k^2}{k!^4} +O(n^{2k-2}) .
\end{eqnarray}

Next we look at ordered pairs of indices $(i_1,i_2)$ so that the corresponding
subwords $I_1$ and $I_2$ intersect in exactly one entry, the entry
 $x$. Let us say that counting
from the left, $x$ is the $a$th 
entry in $I_1$, and the $b$th entry in $I_2$. See Figure
\ref{subwords} for an illustration.

\begin{figure}[ht]
\begin{center}
  \epsfig{file=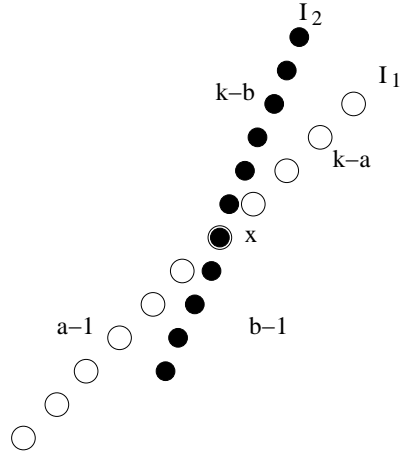}
  \caption{In this example, $k=11$, $a=7$, and $b=5$. }
  \label{subwords}
\end{center}
\end{figure}

Observe that
 $X_{i_1}X_{i_2}=1$ if and only if all of the following independent 
events hold. 
\begin{itemize}
\item In the $(2k-1)$-element set of entries that belong to $I_1\cup I_2$, 
the entry $x$ is the $(a+b-1)$th smallest. This happens with
probability $1/(2k-1)$. 
\item The $a+b-2$ entries on the left of $x$ in  $I_1\cup I_2$ are all smaller
than the $2k-a-b$ entries on the right of $x$ in $I_1\cup I_2$.
This happens with probability $\frac{1}{{2k-2\choose a+b-2}}$. 
\item The subwords of $I_1$ on the left of $x$ and on the right of $x$, 
and the subwords of $I_2$  on the left of $x$ and on the right of $x$
are all monotone increasing. This happens with probability
$\frac{1}{(a-1)!(b-1)!(k-a)!(k-b)!}$. 
\end{itemize}
Therefore, if  $|I_1\cap I_2|=1$, then
\begin{eqnarray} \label{oneprob}
P(X_{i_1}X_{i_2}=1) & = & 
\frac{1}{(2k-1){2k-2\choose a+b-2}(a-1)!(b-1)!(k-a)!(k-b)!} \\
& = & \frac{1}{(2k-1)!}\cdot{a+b-2 \choose a-1}{2k-a-b\choose k-a}
.\end{eqnarray}

How many such ordered pairs $(I_1,I_2)$ are there? There are ${n\choose 2k-1}$
choices for the underlying set $I_1\cup I_2$. Once that choice is made,
the $a+b-1$st smallest entry of $I_1\cup I_2$ will be $x$. Then
the number of choices for the set of  entries other than $x$
that will be part of $I_1$ is ${a+b-2\choose a-1}{2k-a-b\choose k-a}$.
 Therefore, summing over all $a$ and $b$ and
recalling (\ref{oneprob}), 
\begin{eqnarray} \label{contribution}
p_n & = & \sum_{|I_1\cap I_2|=1} E(X_{i_1}X_{i_2}=1) \\
&  = &  
 \frac{1}{(2k-1)!}{n\choose 2k-1}\sum_{a,b}{a+b-2 \choose a-1}^2
{2k-a-b\choose k-a}^2.
\end{eqnarray}
The expression we just obtained is a polynomial of degree $2k-1$, in the
variable $n$. We claim that its leading coefficient is 
larger than $k^2/k!^4$. If we can show that, the proposition will be proved
since (\ref{theeasy}) shows that the summands not included in 
(\ref{contribution}) contribute about $-\frac{k^2}{k!^4}n^{2k-1}$ to 
the left-hand side of (\ref{simplified}). 

Recall that by the Cauchy-Schwarz inequality, if $t_1,t_2,\cdots, t_m$
are non-negative real numbers, then 
\begin{equation}\label{schwarz}
\frac{\left(\sum_{i=1}^m t_i\right)^2}{m} \leq \sum_{i=1}^m t_i^2,
\end{equation}
where equality holds if and only if all the $t_i$ are equal.

Let us apply this inequality with the numbers ${a+b-2 \choose a-1}^2
{2k-a-b\choose k-a}^2$ playing the role of the $t_i$, where $a$ and $b$ 
range from 1 to $k$.
We get that 
\begin{equation}
\label{cauchy} \sum_{1\leq a,b \leq k}{a+b-2 \choose a-1}^2
{2k-a-b\choose k-a}^2  > \frac{\left (
\sum_{1\leq a,b\leq k} {a+b-2 \choose a-1}
{2k-a-b\choose k-a} \right)^2}{k^2}. \end{equation}
We will use Vandermonde's identity to compute the right-hand side. To that
end, we first compute the sum of summands with a {\em fixed} $h=a+b$.
We obtain
\begin{eqnarray}
\sum_{1\leq a,b\leq k} {a+b-2 \choose a-1}
{2k-a-b\choose k-a} &  = & \sum_{h=2}^{2k} \sum_{a=1}^k 
{h-2\choose a-1}{2k-h\choose k-a} \\
 &  = & \sum_{h=2}^{2k} {2k-2\choose k-1} \\
 &  = & (2k-1) \cdot {2k-2\choose k-1}.
\end{eqnarray}
Substituting the last expression into the right-hand side of (\ref{cauchy})
yields
\begin{equation} \label{estimate} \sum_{1\leq a,b \leq k}{a+b-2 \choose a-1}^2
{2k-a-b\choose k-a}^2 >  \frac{1}{k^2} \cdot (2k-1)^2 \cdot
{2k-2\choose k-1}^2.\end{equation}

\noindent Therefore, (\ref{contribution}) and (\ref{estimate}) imply that
\[p_n>\frac{1}{(2k-1)!}{n\choose 2k-1}\frac{(2k-1)^2}{k^2} 
{2k-2\choose k-1}^2.\]

As we pointed out after (\ref{contribution}), $p_n$ is a polynomial of 
degree $2k-1$ in the variable $n$. The last displayed inequality shows that
its leading coefficient is larger than
\[ \frac{1}{(2k-1)!^2} \cdot \frac{1}{k^2} \cdot \frac{(2k-2)!^2}{(k-1)!^4}
=\frac{k^2}{k!^4}  \] as claimed. 

Comparing this with (\ref{theeasy})
 completes the proof of our Proposition.
\end{proof}
 
We can now return to the application of Theorem \ref{janson} to our 
variables $X_{n,i}$. By Proposition \ref{varprop}, there is an absolute
constant $C$ so that $\sigma_n>Cn^{k-0.5}$ for all $n$. 
So (\ref{jansencond}) will be satisfied if we show that
there exists a positive integer $m$ so that
\[{n\choose k} (dn^{k-1})^{m-1} \cdot (n^{-k+0.5})^m<
dn^{-0.5m} 
 \rightarrow 0.\]
Clearly, any positive integer $m$ is a good choice. So we have proved the
following theorem.

\begin{theorem} Let $k$ be a fixed positive integer, 
and let $X_n$ be the random variable counting
 occurrences of $\alpha_k$ in permutations of length $n$. 
Then $\tilde{X}_n\rightarrow N(0,1)$. In other words, $X_n$ is asymptotically
normal.
\end{theorem}

\section{Monotone Subsequences with Entries in Consecutive Positions}
In 2001, Sergi Elizalde and Marc Noy \cite{elizalde}
considered similar problems using another definition of pattern containment.
Let us say that the permutation $p=p_1p_2\cdots p_n$ {\em tightly} 
contains the permutation  $q=q_1q_2\cdots q_k$ if there exists an index
$0\leq i\leq n-k$ so that $q_j<q_r$ if and only if $p_{i+j}<p_{i+r}$. 
(We point out that this definition is a very special case of the one
 introduced
by Babson and Steingrimsson in \cite{babson} and called {\em generalized
pattern avoidance}, but we will not need that much more general concept in 
this paper.) 

\begin{example}
While permutation 246351 contains 132 (take the second, third,
and fifth entries), it does not {\em tightly contain}
132 since there are no three entries in consecutive positions in 246351
that would form a 132-pattern. 
\end{example}

 If $p$ does not tightly contain $q$, then we say that $p$ {\em
tightly avoids} $q$. Let $T_n(q)$ denote the number of $n$-permutations
that tightly avoid $q$. An intriguing conjecture of Elizalde and Noy 
\cite{elizalde} is the following.

\begin{conjecture} \label{enoy}
 For any pattern $q$ of length $k$ and for any positive
integer $n$, the inequality 
\[T_n(q)\leq T_n(\alpha_k)\]
holds.
\end{conjecture}

This is in stark contrast with the situation for traditional patterns, where,
as we have seen in the previous section,
the monotone pattern is not the easier or the harder to avoid, even in the
sense of growth rates. 

\subsection{Tight Patterns of Length Three} 
Conjecture \ref{enoy}
 is proved in \cite{elizalde} in the special case of $k=3$.
As it is clear by taking reverses and complements that $T_n(123)=T_n(321)$ and
that $T_n(132)=T_n(231)=T_n(213)=T_n(312)$, it suffices to show that
$T_n(132)<T_n(123)$ if $n\geq n$. The authors achieve that by a simple
injection. 

It turns out that the numbers $T_n(123)$ are not simply larger than the
numbers $T_n(132)$; they are larger even in the sense of logarithmic
asymptotics. The following results contain the details. 

\begin{theorem} \label{tight123} \cite{elizalde}
 Let $A_{123}(x)=\sum_{n\geq 0} T_n(123)\frac{x^n}{n!}$ be the
exponential generating function of the sequence $\{T_n(123)\}_{n\geq 0}$.
Then
\[A_{123}(x)=\frac{\sqrt{3}}{2} \cdot \frac{e^{x/2}}{\cos \left (
\frac{\sqrt{3}}{2}x+\frac{\pi}{6}\right )}.\] Furthermore,
\[T_n(123)\sim \gamma_1 \cdot (\rho_1)^n \cdot n!,\]
where $\rho_1=\frac{3\sqrt{3}}{2\pi}$ and $\gamma_1=e^{3\sqrt{3}\pi}$.
\end{theorem}

\begin{theorem} \label{tight132} \cite{elizalde}
Let $A_{132}(x)=\sum_{n\geq 0} T_n(132)\frac{x^n}{n!}$ be the
exponential generating function of the sequence $\{T_n(132)\}_{n\geq 0}$.
Then
\[A_{132}(x)=\frac{1}{1-\int_0^x e^{-t^2/2} dt}.\]
Furthermore, 
\[T_n(132)\sim \gamma_2 \cdot (\rho_2)^n \cdot n!,\]
where $\rho_2^{-1}$ is the unique positive root of the equation
$\int_0^x e^{-t^2/2} dt=1$, and $\gamma_2= e^{(\rho_2)^{-2}/2}$.
\end{theorem}

\subsection{Tight Patterns of Length Four}

For tight patterns, the case of length four is even more complex than it is
for traditional patterns in that for tight patterns. Indeed,
 it is not true that
each of the 24 sequences $T_n(q)$, where $q$ is a tight pattern of length 
four, is identical to one of $T_n(1342)$, $T_n(1234)$, and $T_n(1324)$.
In fact, in \cite{elizalde}, Elizalde and Noy showed that there are exactly
seven distinct sequences of this kind. They have also proved the following
results.

\begin{theorem} We have 
\begin{enumerate}
\item $T_n(1342)\sim \gamma_1 (\rho_1)^n \cdot n!$,
\item $T_n(1234)\sim \gamma_2 (\rho_2)^n \cdot n!$, and
\item $T_n(1243)\sim \gamma_3 (\rho_3)^n \cdot n!$,
\end{enumerate}
where $\rho_1^{-1}$ is the smallest positive root $z$ of the equation
$\int_0^z=e^{-t^3/6}dt=1$, $\rho_2^{-1}$ is the smallest positive root
of $\cos z -\sin z +e^{-z}=0$, and $\rho_3$ is the solution of a certain 
equation involving Airy functions. 

The approximate values of these constants are 
\begin{itemize}
\item $\rho_1=0.954611$, $\gamma_1=1.8305194$, 
\item $\rho_2=0.963005$,  $\gamma_2=2.2558142$, 
\item $\rho_3=0.952891$,  $\gamma_3=1.6043282$.
\end{itemize}
\end{theorem}

These results are interesting for several reasons. First, we see that
again, $T_n(\alpha_4)$ is larger than the other $T_n(q)$, even in the
asymptotic sense. Second, $T_n(1234)\neq T_n(1243)$, in contrast to 
the traditional case, where $S_n(1234)=S_n(1243)$. Third, the tight
pattern 1342 is {\em not} the hardest to avoid, unlike in the traditional
case, where $S_n(1342)\leq S_n(q)$ for any pattern $q$ of length four. 

\subsection{Longer Tight Patterns}

For tight patterns that are longer than four, the only known results concern
monotone patterns. They have been found by Richard Warlimont, and, 
independently, also by Sergi Elizalde and Marc Noy. 

\begin{theorem} \label{longer1}
 \cite{elizalde}, \cite{warlimont1}, \cite{warlimont2}
For all  integers $k\geq 3$, the identity
\[\sum_{n\geq 0}T_n(\alpha_k)\frac{x^n}{n!}=\left(\sum_{i\geq 0}
\frac{x^{ik}}{(ik)!} - \sum_{i\geq 0} \frac{x^{ik+1}}{(ik+1)!}
\right )^{-1}
\] holds.
\end{theorem}

\begin{theorem}   \label{longer2} \cite{warlimont2}
Let $k\geq 3$, let $f_k(x)=\sum_{i\geq 0}
\frac{x^{ik}}{(ik)!} - \sum_{i\geq 0} \frac{x^{ik+1}}{(ik+1)!}$, and
let $\omega_k$ denote the smallest positive root of $f_k(x)$.
Then
\[\omega_k=1+\frac{1}{m!}\left(1+O(1)\right),\]
and
\[\frac{T_n(\alpha_k)}{n!} \sim c_m\omega_k^{-n}.\]
\end{theorem}

\subsection{Growth Rates}
The form of the results in Theorems \ref{tight123} and \ref{tight132} is
not an accident. They are special cases of the following general theorem.

\begin{theorem} \cite{sergi} \label{sergi}
For all patterns $q$, there exists a constant $w_q$ so that
\[\lim_{n\rightarrow \infty} \left ( \frac{T_n(q)}{n!} \right )^{1/n}
=w_q.\]
\end{theorem}

Compare this with the result of Corollary \ref{limits}. That Corollary
and the fact that the sequence $(S_n(q)^{1/n}$ is increasing, show
 that the numbers $S_n(q)$ are roughly as large as $L(q)^n$, for some
constant $L(q)$. Clearly, it is much easier to avoid a tight pattern than
a traditional pattern. However, Theorem \ref{sergi} shows how much easier
it is. Indeed, this time it is not the {\em number} of pattern avoiding
permutations is simply exponential; it is their {\em ratio} to all 
permutations that is exponential. 

The fact that $T_n(q)/n! < C_q^n$ for {\em some} $C_q$ is straightforward.
Indeed, $T_n(q)/n! < \left(\frac{k!-1}{k!}\right)^{\lfloor n/k \rfloor}$
by simply looking at $ \lfloor n/k \rfloor$ distinct subwords of $k$ 
consecutive entries. Interestingly, Theorem \ref{sergi} shows that this 
straightforward estimate is optimal in some (weak) sense. Note that there
is no known way to get a result similarly close to the truth for traditional
patterns. 
  
\subsection{Asymptotic Normality} \label{tightass}
Our goal now is to prove that the distribution of tight copies of
 $\alpha_k$ are asymptotically
normal in randomly selected permutations of length $n$. Note that
in the special case of $k=2$, our problem is
reduced to the classic result stating that descents of permutations
are asymptotically normal. (Just as in the previous section, see
\cite{fulman} and its references for various proofs of this fact, or
\cite{ngendes} for a generalization.) Our method is
very similar to the one we used in Subsection \ref{asymptotics}. For
fixed $n$ and $1\leq i\leq n-k+1$, let $Y_{n,i}$
 denote the indicator random variable of the event that in
 $p=p_1p_2\cdots p_n$, the subsequence $p_ip_{i+1}\cdots p_{i+k-1}$ is 
increasing. Set $Y_n=\sum_{i=1}^{n-k+1}Y_{n,i}$.
 We want to use Theorem \ref{janson}.
Clearly, $|Y_{n,i}|\leq 1$ for every $i$, and $N_n=n-k+1$. Furthermore,
the graph with vertex set $\{1,2,\cdots ,n-k+1\}$ in which there is an edge
between $i$ and $j$ if and only if $|i-j|\leq k-1$ is a dependency graph
for the family $\{Y_{n,i}|1\leq i\leq n-k+1\}$. In this graph, 
$\Delta_n=2k-2$. We will prove the following estimate for  $\Var(Y)$.

\begin{proposition} \label{tightupper}
There exists a positive constant $C$ so that $\Var(Y)\geq cn$ for all $n$.
\end{proposition}

\begin{proof}
By linearity of expectation, we have
\begin{eqnarray} \label{tvariance}
\Var (Y_n) & = & E(Y_n^2) - (E(Y_n))^2 \\
 & = & E \left (\left( \sum_{i=1}^{n-k+1} Y_{n,i} \right )^2 \right )
-  \left (E \left (\sum_{i=1}^{n-k+1} Y_{n,i} \right ) \right )^2 \\
 & = & E \left (\left( \sum_{i=1}^{n-k+1} Y_{n,i} \right )^2 \right )
- \left( \sum_{i=1}^{n-k+1} E(Y_{n,i}) \right )^2  \\
 \label{tlastone} & = &  \sum_{i_1, i_2}
E(Y_{n,i_1}Y_{n,i_2})  - \sum_{i_1, i_2}
E(Y_{n,i_1})E(Y_{n,i_2}).
\end{eqnarray} 

In (\ref{tlastone}), all the $(n-k+1)^2$ summands with a negative
sign are equal to $1/k!^2$. Among the summands with a positive sign,
the $(n-2k+1)(n-2k+2)$ summands in which $|i_1-i_2|\geq k$ are equal
to $1/k!^2$, the $n-k+1$ summands in which $i_1=i_2$ are equal to 
$1/k!$, and the $2(n-2k+2)$ summands in which $|i_1-i_2|=k-1$ are 
equal to $1/(k+1)!$. All remaining summands are non-negative. This
shows that
\begin{eqnarray*} \Var(Y_n) & \geq & \frac{n(1-2k)+3k^2-2k+1}{k!^2}
+\frac{n-k+1}{k!} +  \frac{2(n-k+2)}{(k+1)!} \\
& \geq & \left(\frac{1}{k!}+\frac{2}{(k+1)!} - \frac{2k-1}{k!^2}\right) n
+d_k,\end{eqnarray*}
where $d_k$ is a constant that depends only on $k$. As the coefficient
$\frac{1}{k!}+\frac{2}{(k+1)!} - \frac{2k-1}{k!^2}$ of $n$ in the
last expression  is positive for
all $k\geq 2$, our claim is proved. 
\end{proof}

The main theorem of this subsection is now immediate.

\begin{theorem} Let $Y_n$ denote the random variable counting tight copies
of $\alpha_k$ in a randomly selected permutation of length $n$. Then
$\tilde{Y}_n\rightarrow N(0,1)$.
\end{theorem}

\begin{proof} Use Theorem \ref{janson} with $m=3$, and let $C$ be the
constant of Proposition \ref{tightupper}. Then (\ref{jansencond})
simplifies to 
\[(n-k+1)\cdot (2k-2)^2 \cdot \frac{C^3}{n^{1.5}},\] which converges
to 0 as $n$ goes to infinity. 
\end{proof}

\section{Consecutive Entries in Consecutive Positions}
Let us take the idea of Elizalde and Noy
one step further, by restricting the notion of pattern containment further
as follows. Let $p=p_1p_2\cdots p_n$ be a permutation, let $k<n$, and let
 $q=q_1q_2\cdots
q_k$ be another permutation. We say that $p$ {\em very tightly} contains
$q$ if there is an index $0\leq i\leq n-k$ and an integer $0\leq a\leq 
n-k$ so that  $q_j<q_r$ if and only if $p_{i+j}<p_{i+r}$, and, 
\[\{p_{i+1},p_{i+2},\cdots ,p_{i+k}\}=\{a+1,a+2,\cdots ,a+k\} .\]
That is, $p$ very tightly contains $q$ if $p$ tightly contains $q$ and
the entries of $p$ that form a copy of $q$ are not just in consecutive
positions, but they are also consecutive as integers (in the sense that 
their set is an interval). We point out that this definition was used
by A. Myers \cite{myers} who called it {\em rigid} pattern avoidance.
However, in order to keep continuity with our previous definitions, we 
will refer to it as very tight pattern avoidance.

For example, 15324 tightly contains 132 (consider the first three entries),
but does not very tightly contain 132. On the other hand, 15324 very tightly
contains 213, as can be seen by considering the last three entries. If $p$
does not very tightly contain $q$, then we will say that $p$ {\em very tightly
avoids} $q$. 

\subsection{Enumerative Results}

Let $V_n(q)$ be the number of permutations of length $n$ that very tightly
avoid the pattern $q$. The following early results on $V_n(\alpha_k)$ are due
to David Jackson and al. They generalize earlier work by Riordan 
\cite{Riordan} concerning the special case of $k=3$.

\begin{theorem}  \label{jackson} \cite{jackson}, \cite{jackreid}
For all positive integers $n$, and any $k\leq n$, the value of 
$V_n(\alpha_k)$ is equal to 
the coefficient of $x^n$ in the formal power series
\[\sum_{m\geq 0} m!x^m \left(\frac{1-x^{k-1}}{1-x^k}\right)^m.\] 
\end{theorem}

Note that in particular, this implies that for $k\leq n<2k$, the number of
permutations of length $k+r$ {\em containing} a very tight copy of $\alpha_k$
is $r!(r^2+r+1)$.

\subsection{An Extremal Property of the Monotone Pattern}

Recall that we have seen in Section 2 that in the multiset of 
the $k!$ numbers $S_n(q)$ where $q$ is of length $k$,
 the number $S_n(\alpha_k)$ is neither minimal nor maximal. Also recall
that in Section 3 we mentioned that 
in the multiset of 
the $k!$ numbers $T_n(q)$, where $q$ is of length $k$, 
 the number $T_n(\alpha_k)$ is {\em conjectured}
to be maximal. While we cannot prove that  we prove that in the 
in the multiset of 
the $k!$ numbers $V_n(q)$, where $q$ is of length $k$, 
the number $V_n(\alpha_k)$  is maximal, in this Subsection we prove
that for almost all very tight patterns $q$ of length $k$, the inequality
$V_n(q)\leq V_n(\alpha_k)$ does hold.  

\subsubsection{An Argument Using Expectations} \label{outline}
Let $q$ be any pattern of length $k$.
For a fixed positive integer  $n$, let $X_{n,q}$ be the random 
variable counting the very tight copies
of $q$ in a randomly selected 
$n$-permutation. It is straightforward to see that by linearity of 
expectation,
\begin{equation} \label{equal} 
E(X_{n,q})=\frac{(n-k+1)^2}{{n\choose k}k!}.\end{equation}
In particular, $E(X_{n,q})$ does not depend on $q$, just on the length $k$
of $q$.

Let $p_{n,i,q}$
 be the probability that a randomly selected $n$-permutation
contains {\em exactly} $i$ very tight copies of $q$, and let $P(n,i,q)$ be the
probability that a randomly selected $n$-permutation contains {\em at least} 
$i$  very tight 
copies of $q$. Note that $V_n(q)=(1-P(n,1,q))n!$, for any given pattern
$q$. 

Now note that by the definition of expectation
\begin{eqnarray*} \label{atleast} 
E(X_n,q) & = & \sum_{i= 1}^{m} ip_{n,i,q} \\
  & = & \sum_{j=0}^{m-1} \sum_{i=0}^j p_{n,m-i,q} \\
 & = & p_{n,m,q} + (p_{n,m,q}+p_{n,m-1,q}) + \cdots + 
(p_{n,m,q}+\cdots +p_{n,1,q}) \\
 & = & \sum_{i=1}^m P(n,i,q) .
\end{eqnarray*}
We know from (\ref{equal}) that $E(X_{n,q})=E(X_{n,\alpha_k})$,
and then the previous displayed equation implies that
\begin{equation}\label{bigequal}
\sum_{i=1}^m P(n,i,q) = \sum_{i=1}^m P(n,i,\alpha) .\end{equation}
So if we can show that for $i\geq 2$, the inequality 
\begin{equation} \label{toprove} P(n,i,q)\leq P(n,i,\alpha_k)
\end{equation} holds, then (\ref{bigequal}) will imply
that  $P(n,1,q) \geq P(n,1,\alpha_k)$, which is equivalent to $V_n(q)\leq 
V_n(\alpha_k)$, which we set out to prove. 

\subsubsection{Extendible and Non-extendible Patterns}
Now we are going to describe the set of patterns $q$ for which we will
prove that $V_n(q)\leq 
V_n(\alpha_k)$.

Let us assume that the permutation $p=p_1p_2\cdots p_n$ very tightly contains
two {\em non-disjoint} copies of the pattern $q=q_1q_2\cdots q_k$.
Let these two copies be  $q^{(1)}$ and
 $q^{(2)}$, so that $q^{(1)}=p_{i+1}p_{i+2}\cdots p_{i+k}$ and 
$q^{(2)}=p_{i+j+1}p_{i+j+2}\cdots p_{i+j+k}$ for some $j\in [1,k-1]$.
 Then $|q^{(1)}\cap q^{(2)}|=k-j+1=:s$. Furthermore, since the set of entries 
of $q^{(1)}$ is an interval, and the set of entries of $q^{(2)}$ is  an
interval, it follows that the set of entries of $q^{(1)}\cap q^{(2)}$ is
also an interval. So the rightmost
 $s$ entries of $q$, and the leftmost $s$ entries
of $q$ must form identical patterns, and the respective sets of these entries
must both be intervals. 

If $q'$ is the reverse of the pattern $q$, then clearly $V_n(q)=V_n(q')$.
Therefore, we can assume without loss of generality that
 that the first entry of $q$ is less than the last entry of $q$.
For shortness, we will call such patterns {\em rising} patterns. 

 We claim that if $p$ very tightly 
contains two non-disjoint copies  $q^{(1)}$ and
 $q^{(2)}$ of the rising pattern $q$, and $s$ is defined as above, then 
the {\em rightmost} $s$ entries of $q$ must also be the {\em largest} $s$ 
entries of $q$. This can be seen by considering  $q^{(1)}$. Indeed,
 the set of these entries of   $q^{(1)}$ is the 
intersection of two intervals of the same length, and therefore, must
 be an ending segment of
the interval that starts on the left of the other. An analogous argument, 
applied for  $q^{(2)}$, shows that the leftmost $s$ entries of $q$ must
also be the {\em smallest} $s$ entries of $q$. 
So we have proved the following.

\begin{proposition} \label{conditions}
Let $p$ be a permutation that very tightly  contains copies $q^{(1)}$ and
 $q^{(2)}$ of the pattern $q=q_1q_2\cdots q_k$. Let us assume without
loss of generality that
$q$ is rising. Then 
 $q^{(1)}$ and $q^{(2)}$ are disjoint unless all of the following hold.

There exists a positive integer $s\leq k-1$ so that 
\begin{enumerate}
\item the rightmost $s$ entries of $q$ are also the largest $s$ entries 
of $q$,
and the leftmost $s$ entries of $q$ are also the smallest $s$ entries of $q$,
and 
\item the pattern of the leftmost $s$ entries of $q$ is identical to the 
pattern of the rightmost $s$ entries of $q$. 
\end{enumerate}
\end{proposition}

If $q$ satisfies both of these criteria, then two
very tightly contained copies of $q$ in $p$ may indeed intersect. 
For example, the pattern $q=2143$ satisfies both of the above criteria with
$s=2$, and indeed, 214365 very tightly contains two intersecting copies of
$q$, namely 2143 and 4365. 

The following definition is similar to one in \cite{myers}.

\begin{definition} 
Let $q=q_1q_2\cdots q_k$ be a rising pattern 
 that satisfies both conditions of Proposition 
\ref{conditions}
 Then we say that $q$ is {\em extendible}.

If $q$ is rising and not extendible,  then we say that $q$ is non-extendible.
\end{definition}

Note that the notions of extendible and non-extendible patterns are only
defined for rising patterns here. 

\begin{example} The extendible patterns of length four are as follows:
\begin{itemize} \item 1234, 1324 (here $s=1$), 
\item 2143 (here $s=2$). \end{itemize}
\end{example}

Now we are in a position to prove the main result of this Subsection. 

\begin{theorem} \label{almostall} 
Let $q$ be any pattern of length $k$ so that either $q$ or
its reverse $q'$ is non-extendible. Then
for all positive integers $n$, \[V_n(q)\leq V_n(\alpha_k).\]
\end{theorem}

\begin{proof} We have seen in Subsubsection \ref{outline} that it suffices
to prove (\ref{toprove}).

 On the one hand, 
\begin{equation} \label{lowerbound} 
\frac{(n-k-i+2)!}{n!} \leq P(n,i,\alpha_k),\end{equation} since the number 
of $n$-permutations very tightly
containing $i$ copies of
$\alpha$ is at least as large as the number of $n$-permutations very tightly
containing
the pattern $12\cdots (i+k-1)$. The latter is at least as large as the number
of $n$-permutations that very tightly contain
 a $12\cdots (i+k-1)$-pattern in their first $i+k-1$ positions.

On the other hand, 
\begin{equation} \label{upperbound}
  P(n,i,q)\leq {n-i(k-1)\choose i}^2(n-ik)!\frac{1}{n!}.
\end{equation}
This can be proved by noting that if $S$ is the $i$-element set of
starting positions of $i$ (necessarily disjoint)
 very tight copies of $q$ in an $n$-permutation, and $A_S$ is the event 
that in a random permutation $p=p_1\cdots p_n$,
the subsequence $p_jp_{j+1}\cdots p_{j+k-1}$ is a very tight
 $q$-subsequence for
all $j\in S$, then $P(A_S)={n-i(k-1)\choose i}(n-ik)!
\frac{1}{n!}$. The details can be found in \cite{rules}.

Comparing (\ref{lowerbound}) and (\ref{upperbound}), the claim of
the theorem follows.
Again, the reader is invited to consult \cite{rules} for details.
\end{proof}

It is not difficult to show \cite{rules} that the ratio of extendible
permutations of length $k$ among all permutations of length $k$ converges
to 0 as $k$ goes to infinity. So 
Theorem \ref{almostall} covers almost all patterns of length $k$. 

\subsection{The Limiting Distribution of the Number of Very Tight Copies}

In the previous two sections, we have seen that the limiting distribution
of the number of copies of $\alpha_k$, as well as the 
limiting distribution
of the number of tight  copies of $\alpha_k$, is normal. Very tight
copies behave differently.  We will discuss the special case
of $k=2$, that is,  the case of the very tight pattern 12.

\begin{theorem} \label{poisson}
 Let $Z_n$ be the random variable that counts very tight
copies of 12 in a randomly selected permutation of length $n$.
Then $Z_n$ converges a Poisson distribution with parameter $\lambda=1$. 
\end{theorem}

A version of this
 result was proved, in a slightly different setup, by Wolfowitz
 in \cite{wolfowitz}
and by Kaplansky in \cite{kaplansky}.  They used the {\em method of moments},
which is the following.

\begin{lemma} \cite{rucinski}  Let $U$ be a random variable so that
\begin{enumerate} \item for every positive integer $k$, 
the moment $E(U^k)$ exists, and
\item the variable $U$ is completely determined
by its moments, that is, there is no other variable with the same sequence
of moments.
\end{enumerate}
Let $U_1,U_2,\cdots $ be a sequence of random variables, and let us assume
that for all positive integers $k$,
\[\lim_{n\rightarrow \infty} U_n^k =U^k.\]
Then $U_n \rightarrow U$ in distribution. 
\end{lemma}

\begin{proof} (of Theorem \ref{poisson}.)
It is well-known \cite{vonmises}
that the Poisson distribution (with any parameter)
 is determined by its
moments, so the method of moments can be applied to prove convergence
to a Poisson distribution. Let $Z_{n,i}$ be the indicator random variable of
the event that in a randomly selected $n$-permutation $p=p_1p_2\cdots p_n$,
the inequality $p_{i}+1=p_{i+1}$. Then $E(Z_{n,i})=1/n$, and the probability
that $p$ has a very tight copy of $\alpha_k$ for $k>2$ is $O(1/n)$. Therefore,
 we have
\begin{equation} \label{longone} 
\lim_{n\rightarrow \infty} E(Z_n^j)=\lim_{n\rightarrow \infty} 
E\left(\left (\sum_{i=1}^{n-1} Z_{n,i}\right )^j\right)
=\lim_{n\rightarrow \infty} 
E\left(\left (\sum_{i=1}^{n-1} V_{n,i}\right )^j\right),\end{equation}
  where the $V_{n,i}$ are {\em independent}
random variables and each of them takes value 0 with probability $(n-1)/n$, 
and value 1 with probability $1/n$. (See 
 \cite{wolfowitz} for more details.) 
The rightmost limit in the above displayed equation
 is not difficult to compute. Let $t$ be a fixed
non-negative integer. Then the probability that exactly
$t$ variables $V_{n,i}$ take value 1 is ${n-1\choose t}n^{-t}
(\frac{n-1}{n})^{n-t} \sim \frac{e^{-1}}{t!}$. Once we know the $t$-element
set of the $V_{n,i}$ that take value 1,  each of the $t^j$
strings of length $j$ formed from those $t$ variables contributes 1 to
$E(V^j)$. Summing over all $t$, this proves that
\[\lim_{n\rightarrow \infty} 
E\left(\left (\sum_{i=1}^{n-1} V_{n,i}\right )^j\right) =e^{-1}\sum_{t\geq 0}
\frac{t^j}{j!}.
\]
On the other hand, it is well-known that $e^{-1}\sum_{t\geq 1}
\frac{t^j}{j!}$, the $j$th Bell number, is also the $j$th moment of the
Poisson distribution with parameter 1. Comparing this to (\ref{longone}),
we see that the sequence $E(Z_n^j)$ converges to the $j$th moment of
the Poisson distribution with parameter 1. Therefore, by the method of
moments, our claim is proved.
\end{proof}


\begin{thebibliography}{99}
\bibitem{albert} M. Albert, M. Elder, A. Rechnitzer, P. Westcott, 
and M. Zabrocki, A lower bound on the growth rate of the class of
4231 avoiding permutations, {\em Adv. Appl. Math, to appear.}
\bibitem{Arratia} R. Arratia, On the Stanley-Wilf conjecture for the number 
of Permutations avoiding a given pattern. {\em Electronic
 J. Combin.}, {\bf 6} (1999), no. 1,  N1.
\bibitem{babson} E. Babson, E. Steingrgrimsson,
Generalized permutation patterns and a classification of the
Mahonian statistics, {\em Sieminaire Lotharingien de Combinatoire}
{\bf 44} (2000), Article B44b.
\bibitem{Deift}
J. Baik, P. Deift, K. Johansson, 
 On the distribution of the length of the longest increasing
 subsequence of random permutations. 
 {\em  J. Amer. Math. Soc.,} {\bf 12} (1999), no. 4, 1119-1178.
\bibitem{Bona1342}
M. B\'ona, Exact enumeration of 1342-avoiding permutations;
A close link with labeled trees and planar maps.
 {\em J. Combin. Theory A},  {\bf 80} (1997), 257--272.
\bibitem{bonathesis}  M. B\'ona, Exact and Asymptotic Enumeration of 
Permutations with Subsequence Conditions, {\em Ph. D. thesis}, 
Massachusetts Institute of Technology, 1997.
\bibitem{bona} M. B\'ona, {\bf Combinatorics of Permutations}, CRC Press, 
2004.
\bibitem{bonal} M. B\'ona,  The limit of a Stanley-Wilf sequence is
not always rational, and layered patterns beat monotone patterns, {\em
 J. Combin. Theory  Ser. A}  {\bf 110} (2005),  no. 2, 223--235.
\bibitem{records}  M. B\'ona,   New Records on Stanley-Wilf Limits.  {\em 
Europ. J.  
Combin.} {\bf 28} (2007), vol. 1, 75-85.
\bibitem{rules} M. B\'ona, Where the monotone pattern (mostly) rules, 
{\em Discrete Mathematics}, to appear. 
\bibitem{ngendes} M. B\'ona, Generalized Descents and Normality, {\em 
submitted}.
\bibitem{elizalde} S. Elizalde; M. Noy, 
Consecutive patterns in permutations, Formal
 power series and algebraic combinatorics (Scottsdale, AZ, 2001), 
{\it  Adv. in Appl. Math.} {\bf   30}  (2003),  no. 1-2, 110--125.
\bibitem{sergi} S. Elizalde, 
Asymptotic Enumeration of Permutations Avoiding Generalized 
Patterns, {\em Adv. Appl. Math.} {\bf 36} (2006), 138-155.
\bibitem{fulman} J. Fulman, Stein's Method and Non-reversible Markov Chains.
 Stein's method: expository lectures and applications,  69--77, IMS Lecture
 Notes Monogr. Ser., 46, Inst. Math. Statist., Beachwood, OH, 2004.
\bibitem{GesselP}
I. Gessel, Personal communication, 1997.
\bibitem{GesselF}
I. Gessel, Symmetric functions and P-recursiveness.
{\em  J. Combin. Theory Ser. A }, {\bf 53} (1990), no. 2,  257--285. 
\bibitem{jackreid} D. M. Jackson; R. C. Reid,   A note on permutations 
without runs of given length, {\em Aequationes Math}, {\bf 17}
 (1978), no. 2-3, 336-343.
\bibitem{jackson}  D. M. Jackson; J. W. Reilly, Permutations with a 
prescribed number of p-runs, {\em Ars Combinatoria}
 {\bf 1} (1976), no. 1, 297-305.
\bibitem{janson} S. Janson, Normal convergence by higher semi-invariants with 
applications to sums of dependent random variables and random graphs.
{\em Ann. Prob.} {\bf 16} (1988), no. 1, 305-312.
\bibitem{kaplansky} I. Kaplansky, The asymptotic distribution of runs of
 consecutive elements,
{\em Ann. Math. Statistics} {\bf 16} (1945), 200--203.
\bibitem{Marcus}
A. Marcus, G. Tardos,
\newblock Excluded Permutation Matrices and the Stanley-Wilf conjecture.
\newblock  {\em J. Combin. Theory Ser. A}, {\bf 107}  (2004),  no. 1,
 153--160.
\bibitem{myers} A. Myers, Counting Permutations by their Rigid Patterns, 
 {\em J. Combin. Theory  Ser. A}  {\bf 99} (2002)  no. 2, 345--357.
\bibitem{Regev} 
A. Regev, Asymptotic values for degrees associated with
strips of Young diagrams, 
{\em Advances in Mathematics},  {\bf 41} (1981), 115--136.
\bibitem{Riordan} J. Riordan, Permutations without 3-sequences, {\it 
Bull. Amer. Math. Soc.}, {\bf 51} (1945), 745-748.
\bibitem{rucinski} A. Rucinski, Proving Normality in Combinatorics,
in {\em Random Graphs}, Volume 2, Wiley Interscience, 1992, 215--231.
\bibitem{Simion} 
R. Simion,  F. W. Schmidt, Restricted permutations.
 {\em European Journal of Combinatorics}, {\bf 6} (1985), 383--406.
\bibitem{vonmises} R. von Mises, \"Uber die Wahrscheinlichkeit seltener 
Ereignisse, {\em Z. Angew. Math. Mech.} {\bf 1} (1921), 121--124.
\bibitem{warlimont1}  R. Warlimont, Permutations avoiding consecutive
 patterns. II.  {\em Arch. Math.} (Basel)  {\bf 84}
  (2005),  no. 6, 496--502. 
\bibitem{warlimont2} R. Warlimont, 
Permutations avoiding consecutive patterns.  {\em Ann. Univ. Sci. Budapest.
 Sect. Comput.}  {\bf 22}  (2003), 373--393. 
\bibitem{West} 
J. West, {\em Permutations with forbidden subsequences; and,
Stack sortable permutations.} 
 {\em PHD-thesis,} Massachusetts Institute of Technology, 1990.
\bibitem{wolfowitz} J. Wolfowitz, Note on Runs of Consecutive Elements,
Annals Math. Statistics {\bf 15} (1944), 97--98.
\end{thebibliography}
\end{document}